\documentclass[lefteqn,a4paper,11pt]{article}

\usepackage[latin1]{inputenc}
\usepackage[T1]{fontenc}
\usepackage{amsthm,amssymb,amsfonts,latexsym,amsmath,amstext}

\usepackage{euscript}

\usepackage{xspace,times,mathptm}

\catcode`\@=11
\newcommand{\languefrancaise}{
\catcode `\?=\active
\def?{\relax\ifhmode\ifdim\lastskip>\z@\unskip\fi
\kern.2em\fi \string?}

\catcode `\;=\active
\def;{\relax\ifhmode\ifdim\lastskip>\z@\unskip\fi
\kern.2em\fi \string;}

\catcode `\:=\active
\def:{\relax\ifhmode\ifdim\lastskip>\z@\unskip\fi
\penalty\@M\ \fi \string:}

\catcode `\!=\active
\def!{\relax\ifhmode\ifdim\lastskip>\z@\unskip\fi
\kern.2em\fi \string!}

\frenchspacing } \catcode`\@=12

\long\def\xcom#1{}

\newcommand{\ladate}{\space\the\day\ \ifcase\month\or janvier\or f\'evrier
\or mars\or avril\or mai\or juin\or juillet\or aout\or septembre
\or octobre\or novembre\or d\'ecembre\fi \ {\oldstyle\the\year}}

\def\alphabet#1{\ifcase#1\or a\or b\or c\or d\or e\or f\or g\or
h\or i\or j\or k\or l\or m\or n\or o\or p\or q\or r\or s\or t\or
u\or v\or w\or x\or y\or z\fi}

\newcommand{\R}{\mathbb{R}}
\newcommand{\PP}{\mathbb{P}}

\newcommand{\N}{\mathbb{N}}

\newcommand{\Z}{\mathbb{Z}}

\newcommand{\unsur}[1]{{\frac{1}{#1}}}
\def\un#1{{{\,\mathbf{1}}_{({#1})}}}
\def\valabs#1{{\left\vert {#1} \right \vert}}

\def\etp#1{{\left ( {#1} \right )}}
\def\etc#1{{\left [ {#1} \right ]}}

\def\norme#1{{\left \Vert #1 \right \Vert}}
\def\ens#1{{\left\{#1\right\}}}

\def\dessus#1#2{\mathord{\mathop{\kern 0pt #2}\limits^#1}}

\newcommand{\bit}{\begin{itemize}}
\newcommand{\eit}{\end{itemize}}
\newcommand{\ben}{\begin{enumerate}}
\newcommand{\een}{\end{enumerate}}

\newcounter{moncompteur}

\newenvironment{myenumerate}%
{\begin{list}{\arabic{moncompteur}. }{\usecounter{moncompteur}%
\setlength{\leftmargin}{0pt}%
\setlength{\labelwidth}{0pt}%
\setlength{\listparindent}{0pt}%
\setlength{\labelsep}{0pt}}}%
{\end{list}}

\def\bmen{\begin{myenumerate}}
\def\emen{\end{myenumerate}}

\newcommand{\undemi}{\frac{1}{2}}

\def\l{\lambda}

\newcommand{\Frond}{{\mathcal F}}

\newcommand{\ABS}[1]{{\left| #1 \right|}} 
\newcommand{\PAR}[1]{{\left(#1\right)}} 
\newcommand{\SBRA}[1]{{\left[#1\right]}} 
\newcommand{\BRA}[1]{{\left\{#1\right\}}} 
\newcommand{\NRM}[1]{{\Vert #1\Vert}} 
\newcommand{\ind}{\mathrm{1}\hskip -3.2pt \mathrm{I}}

\newcommand{\moyn}[1]{{\left<#1\right>^{(n)}}}

\newcommand{\moybis}[1]{{\left<#1\right>}}
\newcommand{\moynbis}[1]{{\left<#1\right>_2^{(n)}}}

\newcommand{\znb}{Z_n(\be)}

\newtheorem{ethm}{Theorem}

\newtheorem{eprop}[ethm]{Proposition}

\newtheorem{elem}[ethm]{Lemma}

\newtheorem{erem}[ethm]{Remark}

\newcommand{\be}{\beta}

\newcommand{\ep}{\epsilon}

\newcommand{\la}{\lambda}

\newcommand{\te}{\theta}

\newcommand{\Si}{\Sigma}

\newcommand{\Om}{\Omega}
\newcommand{\si}{\sigma}

\newcommand{\dE}{\mathbb{E}}

\newcommand{\dN}{\mathbb{N}}

\newcommand{\dR}{\mathbb{R}}

\newcommand{\dZ}{\mathbb{Z}}

\newcommand{\cG}{{\mathcal G }}

\newcommand{\proofend}{\hfill $\Box{~}$}

\newenvironment{eproof}
               {\noindent {\textbf{Proof.}}}
               {\proofend\\}

\def\r{{\mathbb R}}
\def\e{{\mathbb E}}
\def\p{{\mathbb P}}

\def\n{{\mathbb N}}

\def\df{\, {\buildrel{\rm def} \over =}\, }

\def\no{\noindent}

\def\as{\rm a.s.}
\def\qed{\hfill$\Box$}

\newcommand{\gesp}[1]{\mathbf{E}\etc{#1}}

\newcommand{\espxo}[1]{\mathbb{E}_{x_0}\etc{#1}}
\newcommand{\txo}{\tau_{x_0}}

\begin{document}

\title{Strong disorder for a certain class of directed polymers in a random environment}

\author{Philippe {\sc Carmona}\footnote{ P.Carmona: Laboratoire Jean Leray,
    UMR 6629, Université de Nantes, 92208, F-44322, Nantes cedex 03, e-mail:
philippe.carmona@math.univ-nantes.fr}, Francesco {\sc
Guerra}\footnote{F.Guerra: Dipartimento di Fisica, Universit\`a di
Roma ``La Sapienza'', Instituto Nazionale di Fisica Nucleare,
Sezione di Roma 1, Piazzale Aldo Moro, 2 I-00185 Roma, Italy,
e-mail: francesco.guerra@roma1.infn.it},
 Yueyun {\sc Hu}\footnote{Y.Hu: Laboratoire de Probabilités et Modèles
   Aléatoires (CNRS UMR-7599), Université Paris VI, 4 Place Jussieu, F-75252 Paris cedex 05,
e-mail: hu@proba.jussieu.fr}, and Olivier
\sc{Mejane}\footnote{O.Mejane: Laboratoire de Statistique et
Probabilités, Université Paul Sabatier, 118 route de Narbonne
F-31062 Toulouse cedex 04, France,
e-mail:Olivier.Mejane@lsp.ups-tlse.fr}}

\maketitle

\noindent
\begin{tabular}{p{\linewidth}}
\hline
\begin{abstract}
We study  a model of directed polymers in a random environment
with a positive recurrent Markov chain, taking values in a
countable space $\Si$. The random environment is a family
$(g(i,x), i \geq 1, x \in \Si)$ of independent and identically
distributed real-valued variables. The asymptotic behaviour of the
normalized partition function is characterized: when the common
law of the $g(.,.)$ is infinitely divisible and the Markov chain
is exponentially recurrent   we prove that the  normalized
partition function converges exponentially fast towards zero at
all temperatures.
\end{abstract}\\
\hline
\end{tabular}

\section{Introduction}

In the  model of directed polymers in random environment, we study
a  random Gibbs measure defined on the set of paths (of given
length $n$) of a stochastic process. Usually one choose for the
underlying process a simple random walk on $\dZ^d$ (see for
instance \cite{imbrie-spencer},\cite{sinai} or \cite{comets}) or
$\dR^d$ (see \cite{petermann}). In this paper:
\begin{itemize}
\item The stochastic process is an irreducible Markov chain
$(S_n)_{n\in\n}$ with countable state space $\Sigma$, defined on a
probability space $(\Omega, \Frond,\PP_x, x\in\Sigma)$ with
$\PP_x\etp{S_0=x}=1$. \item The environment is a family
$(g(i,x),i\ge 1, x\in\Sigma)$ of   non-degenerate  i.i.d.  random
variables, distributed as a fixed random variable $g$, defined on
a probability space
$({\Omega}^{(g)},\mathbf{\Frond}^{(g)},\mathbf{P})$, having some
exponential moments
  \begin{equation}
    \label{eq:1}
    \exists \, \beta_0\in(0,+\infty]\,, \forall\,  \valabs{\beta}<\beta_0 :  \qquad \gesp{e^{\beta g}} =e^{\lambda(\beta)} < +\infty\,.
  \end{equation}
\item The random energy is the Hamiltonian, defined on the space
$\Omega_n$ of paths of length $n$ by
$$ H_n(g,\gamma) = \sum_{i=1}^n g(i,S_i)\,.$$
(If $\Pi(\cdot,\cdot)$ is the transition matrix of the chain $S$
then
$$\Omega_n=\ens{\gamma=(\gamma(1),\ldots,\gamma(n)) : \Pi(\gamma(i-1),\gamma(i))>0\,, 2\le i\le n}\,).$$
\item For a given inverse temperature $\be>0$, we introduce the
Gibbs measure $\moyn{.}$ on $\Om_n$ and the normalized partition
function $Z_n(\be)$ according to the definitions:
\begin{eqnarray*}
\moyn{f} &\df& \frac{\dE_{x_0} \PAR{ f(S)e^{\be H_n(g,S)-n\la(\be)}}}{\znb} \, , \\
\znb&\df&\dE_{x_0} \PAR{e^{\be H_n(S)-n\la(\be)}} \, ,
\end{eqnarray*}
\noindent for any bounded function $f$ from $\Om_n$ to $\dR$. We
will denote by $\moynbis{.}$ the product probability measure
$\moyn{.}\otimes\moyn{.}$ on $\Om_n^2$.
\end{itemize}

It is elementary to check that $(\znb)_{n \ge 0}$ is a
$((\cG_n)_{n \ge 0},P)$ positive martingale, if $(\cG_n)_{n \geq
0}$ denotes the natural filtration: $\cG_n=\si \PAR{ g(k,x), 1
\leq k \leq n, x \in \Si }$ for $n \geq 1$ and
$\cG_0=\BRA{\emptyset, \Om^g}$. Hence $\znb \xrightarrow[n\to
\infty]{} Z_{\infty}(\be) \ge 0$ almost surely.

Using the terminology of Comets and Yoshida~\cite{comets-continu},
we say there is {\em weak disorder} if $Z_\infty(\beta)>0$ \as,
and {\em strong disorder} if \as\, $Z_\infty(\beta)=0$.

When $(S_n)$ is the simple random walk on $\Z^d$ and when the
environment $g$ is Gaussian, the picture is the following :
\begin{itemize}
\item if $d\ge 3$ and $\beta>\beta_1$ for some $\beta_1>0$, there
is strong disorder and  almost surely $Z_n(\beta)$ converges to
zero exponentially fast.
 \item if $d\ge 3$ and
$\beta<\beta_2$ for some    $\beta_2>0$, then there is weak
disorder. \item if $d=1,2$ then for any $\beta>0$ there is strong
disorder (see \cite{carmona-hu,comets})  with exponential
convergence of $\znb$ to $0$ if $\beta$ is large enough,   but the
rate of convergence is still unknown for small $\beta$.
\end{itemize}

It is not difficult to prove, by the method of  second moment,
that there is weak disorder for a ``transient'' Markov chain when
$\beta$ is small, here by ``transient'' we mean   that $\sum_{n,x}
\PP_{x_0}\etp{S_n=x}^2 < +\infty$.

The aim of this paper is to prove that for a large class of
positive recurrent Markov chain, and for fairly general random
environments,  almost surely $Z_n(\beta)$ converges to zero
exponentially fast.

>From now on, we shall assume that the Markov chain
$(S_n)_{n\in\N}$ is positive recurrent, and that the first return
time to $x_0$,  $\tau_{x_0}=\inf\ens{n\ge 1 : S_n=x_0}$, has small
exponential moments
\begin{equation}
  \label{eq:em}
  \tag{EM} \exists x_0\in\Sigma\,,\exists \kappa>0\,,\quad \espxo{e^{\kappa \txo}}<\infty
\end{equation}

Define $$  p_n(\beta) = \unsur{n}\gesp{\log(\znb)}, \qquad 0\le
\beta < \beta_0, \qquad \beta_0 \in (0, +\infty].$$

Our main result is the following theorem:

\begin{ethm} \label{countable_case}
 If the Markov chain
$(S_n)_{n \ge 0}$ is irreducible, positive recurrent and satisfies
\eqref{eq:em} and if the law of the random environment is
infinitely divisible and satisfies \eqref{eq:1}, then
\begin{itemize}
    \item[(a)]  for small $\beta>0$, the free energy $p(\beta) =
    \lim_{n\to\infty } p_n(\beta)$ exists and $$ {1\over n}  \log
    Z_n(\beta)  \, \to \, p(\beta), \qquad \mbox{a.s. and in
    $L^1$}.$$
    \item[(b)]  the
function   $\beta \in [0, \beta_0)  \to p_n(\beta) $ is non
increasing.
    \item[(c)]  for all  $0< \beta < \beta_0$, $$\limsup_{n\to\infty} p_n(\be)<0\,, \quad \forall\, \beta \in(0,\beta_0) \, .$$
In particular, for any $0<\beta<\beta_0$,  almost surely
$Z_n(\beta)$ converges to zero exponentially fast.
\end{itemize}
\end{ethm}

 This paper is inspired by the works of
Francesco Guerra and Fabio Toninelli (see
\cite{guerra-toninelli}), who developed an interpolation technique
to study the high temperature behaviour of the
Sherrington-Kirkpatrick mean field spin glass model. The principal
ingredient of the proof on the exponential decay  is the
interpolation between the random Hamiltonian $H_n(g,\gamma)$ and a
deterministic Hamiltonian.

\medskip
The paper is organized as follows:
\begin{itemize}
\item In Section 2,  we evaluate the exponential moments of some
additive funstionals whose first consequence is the existence of
the free energy $p(\beta)$ for small $\beta
>0$.   Concentration of measure implies then the \as\, convergence
$\unsur{n}\log \znb \to p(\beta)$.
 \item We devote   Section 3 to
an integration by parts formula, a feature of infinitely divisible
distributions, which entails the monotonicity of free energy (b).

\item The last section contains the proof of Theorem 1.

\end{itemize}

Unless stated otherwise, we assume in the sequel  that $\beta \in
[0, \beta_0)$ and the random environment $g$ is centered.

\section{Exponential moments} \label{lemme_GD}
Recall that  $(S_n)$ is a Markov chain taking values in a
countable set $\Sigma$ satisfying~\eqref{eq:em}, and the
environment variables $(g(i,x))$ are centered and have some
exponential moments (see \eqref{eq:1}).

Let us omit the dependence on  $x_0$ of $\tau$ and denote  the
successive return times to $x_0$ by $\tau_0=0< \tau_1< \tau_2<...<
\tau_n<...$ For a bounded function $f:\Sigma\to\r$, we define
$\overline f = \sup_{x\in\Sigma} f(x)$, $\underline f= \inf_{x\in
\Sigma} f(x)$ and $\| f\|_\infty= \sup_{x\in \Sigma} |f(x)|$. The
main result of this section is the following theorem:

\begin{ethm} \label{T: el}  Let $f: \Sigma \to \r$ be a bounded
measurable function and $ |\beta|< \beta_0 $ is sufficiently small
such that $   \lambda(\beta)+ 2\| f\|_\infty < \kappa.$
\\ (i)  There exists a unique real
number $c(\beta, f) \in  [\underline f,  \overline f +
\lambda(\beta)] $ such that
$$ \unsur{n} \log \e_{x_0} \exp\Big( \beta \sum_{i=1}^{\tau_n} g(i,
S_i) + \sum_{i=1}^{\tau_n} f(S_i) - c(\beta, f) \, \tau_n\Big) \,
\to\,
0, \qquad \mbox{a.s. and in $L^1$}.$$ \\
(ii) We have
$$ \lim_{n\to\infty} {1\over n} \log \e_{x_0} e^{\beta \sum_{i=1}^n  g(i, S_i) +
\sum_{i=1}^n f(S_i)} =  c(\beta,f), \qquad \mbox{a.s. and in
$L^1$}.$$
\end{ethm}

\bigskip
The constant $c(\beta, f)$ does not depend on the starting point
$x_0$, see the forthcoming Remark \ref{R:9}. Taking $\beta=0$ in
Theorem \ref{T: el}, we can evaluate the following Varadhan's type
integral

\begin{eprop}\label{pro:large-devi-result}  For any bounded function $f: \Sigma
\to \r$ such that $ \| f\|_\infty < \kappa/2$, we have $$
\lim_{n\to\infty} {1\over n} \log \e_{x_0} e^{ \sum_{i=0}^{n-1}
f(S_i)} =   c(f),$$ where $c(f) \in [\underline f,    \overline f]
$ is the unique real number such that
$$ \e_{x_0} \exp\Big( \sum_{i=0}^{\tau(x_0)-1} f(S_i) - c(f) \tau(x_0)
\Big) =1.$$
\end{eprop}

According to the theory of large deviations, Proposition
\ref{pro:large-devi-result} is well-known at least for the case
when $(S_n)$ is a Markov chain with finite states, for example by
combining
 Dembo and Zeitouni (\cite{DZ98}, pp. 75) and Ney and Nummelin (\cite{NN87}, Lemma 4.1).  See also  de
Acosta and Ney (\cite{DeA98}) and the references therein for the
large deviation principles  for a Markov chain.

Taking $f=0$ in Theorem \ref{T: el}, we obtain  the existence of
the free energy at high temperature (recalling that $g$ is
centered):

\begin{eprop} \label{P:p4}  Let $ |\beta|< \beta_0 $ be sufficiently small such  that
$  \lambda(\beta) < \kappa.$
\\ (i)  There exists a unique real
number $c(\beta) \in [0, \lambda(\beta)] $ such that $$ {1\over n}
\log \e_{x_0} \exp\Big( \beta \sum_{i=1}^{\tau_n} g(i, S_i) -
c(\beta) \, \tau_n\Big) \, \to\,
0, \qquad \mbox{a.s. and in $L^1$}.$$ \\
(ii) We have
$$ \lim_{n\to\infty} {1\over n} \log \e_{x_0} e^{\beta \sum_{i=1}^n  g(i, S_i)  } =  c(\beta), \qquad \mbox{a.s. and in
$L^1$}.$$
\end{eprop}

Before entering into the proof of Theorem \ref{T: el}, we
establish  a preliminary  result on the concentration of measure,
which  is essentially adapted from Comets, Shiga and Yoshida
(\cite{comets}, Proposition 2.9).

\begin{elem}{\bf (Concentration of measure)}\label{L:cm}  Let $f: \Sigma \to \r$ be a bounded
measurable function and $ |\beta|< \beta_0 $.    Denote by
$D_n(\beta,f)=\beta \sum_{i=1}^n g(i, S_i)+ \sum_{i=1}^n f(S_i)$ and $\overline f =\sup_{x\in \Sigma} f(x)$. \\
(i) Assume that $\overline f + \lambda(\beta) < \kappa$. For any
$\epsilon>0$, there exists a $n_1=n_1(\beta, f, \epsilon)<\infty$
such that for all $n\ge n_1$,
$$ {\bf P}\left( \Big| {1\over n} \log \e_{x_0} e^{ D_{\tau_n}} - {1\over n}{\bf E}\log \e_{x_0} e^{ D_{\tau_n}}  \Big|> \epsilon\right) \le e^{- \epsilon^{2/3} n^{2/3}/6}.$$ \\
(ii)  For any $\epsilon>0$, there exists a $n_2=n_2(\beta, f,
\epsilon)<\infty$ such that for all $n\ge n_2$,
$$ {\bf P}\left( \Big| {1\over n} \log \e_{x_0} e^{ D_n} - {1\over n}{\bf E}\log \e_{x_0} e^{ D_n}  \Big|> \epsilon\right) \le e^{- \epsilon^{2/3} n^{2/3}/4}.$$ \\
(iii)  Assume that $\overline f + \lambda(\beta) < \kappa$ and fix
$1\le a \le b< \infty$.   Then for any $\epsilon>0$, there exists
a $n_3=n_3(a, b, \beta, f, \epsilon)<\infty$ such that for all
$n\ge n_3$,
$$ {\bf P}\left( \Big| {1\over n} \log \e_{x_0} \Big( e^{ D_{\tau_n}} \, \big| \, a n \le \tau_n \le  bn \Big) - {1\over n}{\bf E}\log \e_{x_0} \Big( e^{ D_{\tau_n}}\, \big| \, a n \le \tau_n \le  bn \Big)  \Big|> \epsilon\right) \le e^{- \epsilon^{2/3} n^{2/3}/5},$$
with convention $\e_{x_0}\big(\cdot \big| \emptyset\big)\equiv 1$.
\end{elem}

{\noindent\bf Proof of Lemma \ref{L:cm}:}  Using  the same
 arguments (martingale decomposition, large deviations for martingale)   as that of    Comets, Shiga and Yoshida~\cite{comets}  pp. 720--721,
we obtain (ii) and the following inequality: For any $\epsilon>0$
and $b>1$, there exists a $n_4=n_4(  b, \beta, f,
\epsilon)<\infty$ such that for all $n\ge n_4$ with $\PP_{x_0}(
\tau_n = k(n)) >0$ and $k(n) \le b n$,
\begin{equation} \label{takn} {\bf P}\left( \Big| {1\over n} \log \e_{x_0} \Big( e^{ D_{\tau_n}}1_{(  \tau_n = k(n))}\Big) - {1\over n}{\bf E}\log \e_{x_0} \Big( e^{ D_{\tau_n}}1_{(  \tau_n =k(n))}\Big)  \Big|> \epsilon\right) \le e^{- \epsilon^{2/3}
n^{2/3}/4}.
\end{equation}
Observe that $$ 0\le \log \e_{x_0} \Big( e^{ D_{\tau_n}}1_{( a n
\le \tau_n \le  bn)}\Big) - \max_{an \le k \le b n}  \log \e_{x_0}
\Big( e^{ D_{\tau_n}}1_{(  \tau_n = k)}\Big) \le \log b + \log
n,$$

 \noindent and  for any $u>0$, \begin{eqnarray*} && \left\{
\Big| \max_{an \le k \le b n}  \log \e_{x_0} \Big( e^{
D_{\tau_n}}1_{( \tau_n = k)}\Big) - \max_{an \le k \le b n}{\bf E}
\log  \e_{x_0} \Big( e^{ D_{\tau_n}}1_{( \tau_n = k)}\Big) \Big| >
u \right\}  \\
    && \subset \, \bigcup_{ a n \le k \le  b n}  \left\{
\Big|    \log \e_{x_0} \Big( e^{ D_{\tau_n}}1_{( \tau_n = k)}\Big)
-  {\bf E} \log  \e_{x_0} \Big( e^{ D_{\tau_n}}1_{( \tau_n =
k)}\Big) \Big| > u \right\}
\end{eqnarray*}

The above two observations together with (\ref{takn}) imply (iii).
To prove (i), we remark   that
\begin{equation}\label{bj15} \lim_{b\to\infty} \limsup_{j\to\infty} {1\over j} \log
{\bf E}  \left(\e_{x_0} \Big( e^{ D_{\tau_j}} 1_{(\tau_j \ge b\,
j)}\Big) \right)  = -\infty. \end{equation}

In fact, we have from Fubini's theorem  and Chebychev's inequality
that

\begin{eqnarray*} {\bf E}  \left(\e_{x_0} \Big( e^{ D_{\tau_j}} 1_{(\tau_j \ge b\,
j)}\Big) \right)  &\le & \e_{x_0} \Big( e^{ (\lambda(\beta) +
\overline f   ) \tau_j } 1_{(\tau_j
> b\, j)}\Big) \\
    &\le &  e^{- \delta_0\, b\, j} \Big[ \e_{x_0} e^{ ( \lambda(\beta) + \overline f + \delta_0) \tau}  \Big]^j,
\end{eqnarray*}

\noindent where $\delta_0>0$ denotes a small constant such that $
\lambda(\beta) + \overline f  + \delta_0< \kappa$. This yields
(\ref{bj15}).  Finally,   applying (iii) to $a=1$ (since
$\tau_n\ge n$) and a sufficiently large $b>0$, we obtain (i). \qed

\bigskip

{\noindent\bf Proof of Theorem \ref{T: el}:}  (i) Let $c_+= \kappa
- \lambda(\beta) - \overline f >0$. We shall show  that the
following function $\psi: (-\infty , c_+) \to \r$ is well-defined:
for any $-\infty< c< c_+$,
$$ \psi(c) \df
 \lim_{n\to\infty} (\mbox{a.s. and in $L^1$}){1\over n} \log \e_{x_0} \exp\Big( \beta
\sum_{i=1}^{\tau_n} g(i, S_i) + \sum_{i=1}^{\tau_n} f(S_i) + c  \,
\tau_n\Big) . $$

To this end, we shall apply  the subadditivity theorem. For
notational convenience, denote by $$D_n=D_n(g, S)=   \beta
\sum_{i=1}^n g(i, S_i) + \sum_{i=1}^n f(S_i) + c  \,
 n .$$ Using the strong Markov property at $\tau_n$, we
have $$ \e_{x_0} e^{D_{\tau_{n+m}}} = \sum_j \, \e_{x_0} \Big(
e^{D_{\tau_n}}1_{(\tau_n=j)}\Big)\, \e_{x_0}
e^{D_{\tau_m}(\theta_j g, S)},$$ where $\theta_j$ denotes the
shift operator on $g$: $\theta_j g(i,x)= g(i+j, x)$. By concavity,
\begin{eqnarray*} \log \e_{x_0} e^{D_{\tau_{n+m}}} &=& \log \e_{x_0}
e^{D_{\tau_n}} + \log \sum_j { \e_{x_0} \Big(
e^{D_{\tau_n}}1_{(\tau_n=j)}\Big) \over \e_{x_0} e^{D_{\tau_n}}
}\, \e_{x_0} e^{D_{\tau_m}(\theta_j g, S)} \\
    &\ge& \log \e_{x_0} e^{D_{\tau_n}} + \sum_j {  \e_{x_0}
\Big( e^{D_{\tau_n}}1_{(\tau_n=j)}\Big) \over \e_{x_0}
e^{D_{\tau_n}} }\, \log \e_{x_0} e^{D_{\tau_m}(\theta_j g, S)}.
\end{eqnarray*}

Hence $${\bf E} \log \e_{x_0} e^{D_{\tau_{n+m}} } \ge {\bf E} \log
\e_{x_0} e^{D_{\tau_n} } + {\bf E}\log \e_{x_0} e^{D_{\tau_m}},$$
and $$ \psi(c) \df \lim_{n\to\infty}{1\over n} {\bf E}\log
\e_{x_0} e^{D_{\tau_n} } = \sup_{n\ge1}{1\over n} {\bf E}\log
\e_{x_0} e^{D_{\tau_n} }.$$

The a.s. convergence follows from Lemma \ref{L:cm} (i) since
$\overline f + c + \lambda(\beta) < \kappa$.

As limit of convex and nondecreasing functions, $\psi(\cdot)$ is
convex and nondecreasing. Moreover, $\psi: (-\infty, c_+) \to \r$
is strictly increasing since $\tau_n\ge n$.   By Jensen's
inequality,
$$ \psi(c) \le \log \e_{x_0} e^{ (\lambda(\beta) + c+  \overline f )
\tau_1},$$ which implies that $\psi(- (\overline f +
\lambda(\beta))) \le 0$. Again using Jensen's inequality and the
fact that $g$ is centered, we have
$$ \psi(c) \ge {\bf E} \log \e_{x_0} e^{ D_{\tau_1}} \ge {\bf E}
\log \e_{x_0} e^{ \beta \sum_1^{\tau_1} g(i, S_i) + (c+ \underline
f)\tau_1} \ge (c+ \underline f)\, \e_{x_0} \tau_1,
$$

\noindent hence $\psi(- \underline f)\ge0$. It follows that there
exists a unique real number $c=   c(\beta, f)\in [\underline f,
\overline f + \lambda(\beta)]$ such that $\psi(-c)=0$, proving
(i).

 \medskip

(ii) Define $$ D_n= D_n(g, S) = \beta\sum_{i=1}^n g(i, S_i) +
\sum_{i=1}^n f(S_i) - c(\beta, f) \, n.$$

Then by (i), \begin{equation}\label{tauj} {1\over j} \log \e_{x_0}
e^{ D_{\tau_j}} \, \to\, 0, \qquad \mbox{a.s. and in $L^1$}.
\end{equation}

We are going to prove that
\begin{equation}\label{asdn}\lim_{n\to\infty} {1\over n} \log \e_{x_0}
e^{D_n} =0, \qquad \mbox{a.s.} \end{equation}

It is not difficult to show that the family $ ({1\over n} \log
\e_{x_0} e^{D_n}, n\ge1)$ is bounded in $L^2$, in fact, by
Jensen's inequality, $ {1\over n} \log \e_{x_0} e^{D_n} \ge {\beta
\over n} \e_{x_0} \sum_1^n g(i, S_i) + \underline f - c(\beta,
f).$  On the other hand, since the function $x (\in \r_+) \to
\log^2(x+e)$ is concave,
$$ {\bf E} \Big( \max(0, {1\over n}  \log \e_{x_0} e^{D_n} )\Big)^2
\le {1\over n^2}\log^2 (e+ \, {\bf E} \log \e_{x_0} e^{D_n}) =
O(1).
$$

Therefore, the family $ ({1\over n} \log \e_{x_0} e^{D_n}, n\ge1)$
is  uniformly integrable, which  in view of (\ref{asdn}) implies
that $ \unsur{n} {\bf E} \log \e_{x_0} e^{D_n} \to 0$. This proves
 the $L^1$ convergence part of (ii).

It remains to show (\ref{asdn}) whose proof is divided into two
parts.

{\bf Upper bound  of (\ref{asdn}):} Notice that $\tau_j \ge j$;
therefore, we have
\begin{eqnarray*} \e_{x_0} e^{D_n} &=& \sum_{j=0}^{n-1} \e_{x_0}
\Big(  e^{D_n} 1_{(\tau_j < n\le \tau_{j+1})}\Big) \\
    &=&  \sum_{j=0}^{n-1} \e_{x_0}
\Big(  e^{D_{\tau_j}} 1_{(\tau_j < n)} \e_{x_0} \, \big[
e^{D_k(\theta_{n-k}g,S)} 1_{(k\le \tau )}\big]\big|_{k=
n-\tau_j}\Big) \\
    &\le&  M_n\, \sum_{j=0}^{n-1} \e_{x_0}
\Big(  e^{D_{\tau_j}} \Big),
\end{eqnarray*}

\noindent where $$ M_n = \max_{1\le k\le n} \e_{x_0} \, \big[
e^{D_k(\theta_{n-k}g,S)} 1_{(k\le \tau_1  )}\big].$$  Observe that
\begin{eqnarray*}  {\bf E} M_n &\le &\sum_{k=1}^n {\bf E} \left(\e_{x_0} \, \big[
e^{D_k(\theta_{n-k}g,S)} 1_{(k\le \tau_1  )}\big] \right)  \\
         &\le&
\sum_{k=1}^n \e_{x_0} \big[ e^{ (\lambda(\beta)+  \overline f -
c(\beta, f)) k} 1_{(k \le \tau_1)}\big]  \\
    &\le&  C \,n\, , \end{eqnarray*}
with $$ C=  \e_{x_0} \big[ e^{ (\lambda(\beta)+  \overline  f -
c(\beta, f)) \tau_1 }\big] \le     \e_{x_0} \big[ e^{
(\lambda(\beta)+  2  \| f\|_\infty ) \tau_1 }\big]    < \infty.$$

\noindent By Borel-Cantelli's lemma, almost surely for all large
$n$, $$ M_n \le n^3.$$ This together with the a.s. convergence in
(\ref{tauj}) imply the upper bound: \begin{equation}\label{upbd}
\limsup_{n\to\infty} {1\over n} \log \e_{x_0} e^{D_n} \le 0,
\qquad \mbox{a.s.} \end{equation}

{\bf Lower bound of (\ref{asdn}):}  By means of (\ref{bj15}), for
sufficiently large $b>0$,
$${\bf E}\left( \e_{x_0} \Big( e^{
D_{\tau_j}} 1_{(\tau_j \ge b\, j)}\Big) \right)  \le e^{-2 j}, $$
which     in view of  Borel-Cantelli's lemma  yields that ${\bf
P}(d\omega)$ a.s. for all large $n \ge n_0(\omega)$,

\begin{equation} \label{exx} \e_{x_0} \Big( e^{
D_{\tau_j}} 1_{(\tau_j  \ge b\, j)}\Big) \le e^{-j}.
\end{equation}

 Then by (\ref{tauj}) and (\ref{exx}),
a.s. for all large $j\ge j_0(\epsilon, \omega)$, \begin{equation}
\label{exy} \liminf_{j\to\infty}\, {1\over j} \log \e_{x_0} \Big(
e^{D_{\tau_j}} \, 1_{(\tau_j < b \, j)}\Big) \, \ge\, 0.
\end{equation}

\noindent Let $\epsilon >0$ be small. We divide the interval  $[1,
b]$ into $K=K(\epsilon) = [b/\epsilon]$ intervals $[a_1, a_2),
..., [a_{K-1}, a_K)$ with $a_1=1, a_K=b$ and $a_{k+1} - a_k=
{b-1\over K} < \epsilon$ for $1\le k\le K-1$. For any $0\le k \le
K-1$, we may repeat   the similar  argument of subadditivity  in
(i) and apply    the concentration of measure (Lemma \ref{L:cm},
(iii)). This yields  that
\begin{equation}\label{subad} {1\over j} \log \e_{x_0} \Big(
e^{D_{\tau_j}} \, 1_{(a_k \, j \le \tau_j \, \le a_{k+1} \,
j)}\Big) \, \to\, \gamma_k\, \quad \mbox{a.s. and in $L^1$},
\end{equation}

\noindent for some deterministic constant $\gamma_k \in [- \infty,
0]$ ($\gamma_k \le0$ because of (\ref{tauj})). Note that
$\gamma_k=-\infty$ if and only if for all $j\ge1$,
$\p_{x_0}(\tau_j \in [j a_k, j a_{k+1}])=0$.

We claim that \begin{equation} \label{gamma} \max_{0\le k < K}
\gamma_k=0.\end{equation}

\noindent Otherwise, since $ \gamma_k<0$  for each $k< K$,
$\e_{x_0} \Big( e^{D_{\tau_j}} \, 1_{(a_k \, j \le \tau_j \, \le
a_{k+1} \, j)}\Big) $ converges to $0$ exponentially fast; then
$\e_{x_0} \Big( e^{D_{\tau_j}} \, 1_{(\tau_j < b \, j)}\Big) =
\sum_{k=0}^{K-1} \, \e_{x_0} \Big( e^{D_{\tau_j}} \, 1_{(a_k \, j
\le \tau_j \, \le a_{k+1} \, j)}\Big)$ would also converge to $0$
exponentially fast, which is in contradiction  with (\ref{exy}).
Then we proved (\ref{gamma}).

Now, we proceed to show  the lower bound. Choose a fixed $k \in
[0, K-1]$ such that $\gamma_{k}=0$.  Let $j= [ {n \over
a_{k+1}}]$. We have
\begin{eqnarray} \e_{x_0} e^{D_n} & \ge &   \e_{x_0}
\Big(  e^{D_n} 1_{( j\, a_{k }\,   \le \tau_j < j\, a_{k +1}\,  )}\Big) \nonumber \\
    &=&    \e_{x_0}
\Big(  e^{D_{\tau_j}} 1_{( j\, a_{k }\,   \le \tau_j < j\, a_{k
+1}\,  )} \,  \e_{x_0} \, \big[ e^{D_\ell(\theta_{n-\ell}g,S)}
 \big]\big|_{\ell=n-\tau_j}\Big)  \nonumber\\
    &\ge & m_n \, \e_{x_0}
\Big(  e^{D_{\tau_j}} 1_{( j\, a_{k }\,   \le \tau_j < j\, a_{k
+1}\,  )}\Big), \label{dl}
\end{eqnarray}

\noindent where by our choice of $j$ and $a_k$, $\ell = n- \tau_j
\le n- j a_k \le 2 \epsilon n $ and $$ m_n = \min_{1\le \ell \le 2
\epsilon n} \e_{x_0} \, \big[ e^{D_\ell(\theta_{n-\ell}g,S)}
 \big] . $$

  By Jensen's inequality,
$$ \e_{x_0} \, \big[ e^{D_\ell(\theta_{n-\ell}g,S)}  \big] \, \ge \,   e^{
\e_{x_0} \big[ D_\ell(\theta_{n-\ell}g,S) \,  \big]}.$$

\noindent Since  $  f(x)  - c(\beta, f)  \ge \underline f -
c(\beta, f) \ge -     ( 2 \| f\|_\infty + \lambda(\beta)) > -
\kappa$, we have
$$ \e_{x_0} \big[
D_\ell(\theta_{n-\ell}g,S)  \big] \ge \beta \sum_{i=1}^\ell \sum_x
g(i+n-\ell, x)\, q_i(x) -     \kappa \, \ell  ,
$$

\noindent  where we write  $q_i(x)= \PP_{x_0}\big( S_i= x \big)$
for notational convenience. Observe that $\sum_{i=1}^\ell \sum_x
(q_i(x))^2 = \sum_1^\ell \p_{x_0} \big( S_i= \widetilde
S_i\big)\le \ell \le 2 \epsilon n$, where $\widetilde S$  denotes
an independent copy of $S$. By Chebychev's inequality, for any $v>
0$,
$$ {\bf P} \Big( \sum_{i=1}^\ell \sum_x g(i+n-\ell, x)\, q_i(x ) <
- n^{2/3} \Big) \le e^{ - v\, n^{2/3}}\, e^{ \sum_1^\ell \sum_x
\lambda ( -v\,  q_i(x ))} \le e^{ - {n^{1/3} \over 2
\lambda''(0)}},$$

\noindent where in the last inequality, we choose $v=
{n^{-1/3}\over \lambda''(0)}$ and use the fact that $\lambda (u)
\sim {\lambda''(0)\over2} u^2$ for small $u$. It turns out that $$
{\bf P} \Big( \min_{1\le \ell \le 2 \epsilon n} \sum_{i=1}^\ell
\sum_x g(i+n-\ell, x)\, q_i(x ) < - n^{2/3} \Big) \le n\,  e^{ -
{n^{1/3} \over 2 \lambda''(0)}},$$

\noindent  whose sum on $n$ converges. Hence ${\bf P}$ a.s. for
all large $n\ge n_0(\omega)$, $\min_{1\le \ell \le 2 \epsilon n}
\sum_{i=1}^\ell \sum_x g(i+n-\ell, x)\, q_i(x )  \ge -n^{2/3}$ and
therefore
$$ m_n \ge e^{ - \beta n^{2/3} - 2 \epsilon  \, \kappa \, n} \, .$$

\no Plugging this into (\ref{dl}) and using (\ref{subad}) with
$\gamma_k=0$ by our choice of $k$, we obtain that a.s.
$$ \liminf_{n\to\infty} {1\over n} \log \e_{x_0} e^{D_n} \, \ge -2
\epsilon \, \kappa ,$$ for any $\epsilon>0$.  The lower bound of
(\ref{asdn}) follows by letting $\epsilon \to0$. This together
with the upper bound (\ref{upbd}) complete  the proof of Theorem
\ref{T: el}.
 \qed

\begin{erem}  When $\beta=0$,   the value  of   $a_k$ in (\ref{dl}) can be  easily determined  by  a change of probability
measure.
\end{erem}

\bigskip

\noindent  We shall need the following corollary:
\begin{elem} \label{travail_Carmona}
 Assume $(EM)$.
Let $f$ be a bounded function from $\Si$ to $\dR$. Then  for all
$\ABS{t}<t_0=\frac{\kappa}{2 \norme{f}_\infty}$,  the limit
$$c(t)=\lim_{ n \to \infty }\frac{1}{n} \log \dE_{x_0}
\PAR{e^{t\sum_{i=1}^n f(S_i)}}$$ \noindent exists. Moreover, $c$
is differentiable at $0$ with $c'(0) = \sum_{ x \in \Si} f(x)
\mu(x)$, where $\{\mu(x), x\in \Sigma\}$ denotes the invariant
probability measure of $S$.
\end{elem}
\begin{eproof}

Indeed,  for $\valabs{t}<t_0$, the limit $c(t)=\lim_{ n \to \infty
}\frac{1}{n} \log \dE_{x_0} \PAR{e^{t\sum_{i=1}^n f(S_i)}}$
exists. It is the unique real number $c$ such that $\phi(c,t)=1$
where $\phi$ is the function
$$\phi(c,t)= \e_{x_0} \exp\Big( t \sum_0^{\tau(x_0)-1} f(S_i) - c \,  \tau(x_0)
\Big)\,.$$ Since $\phi$ is continuously differentiable in
$(-t_0,t_0)\times J$ with $J$ and open interval, with derivatives
\begin{gather*}
\frac{\partial \phi}{\partial c}=-\e_{x_0}\etc{ \tau(x_0)
\exp\Big( t \sum_0^{\tau(x_0)-1} f(S_i) - c(f) \tau(x_0) \Big)}
\\
\frac{\partial \phi}{\partial t}=\e_{x_0}\etc{
\sum_0^{\tau(x_0)-1} f(S_i) \exp\Big( t \sum_0^{\tau(x_0)-1}
f(S_i) - c(f) \tau(x_0) \Big)}
\end{gather*}
the implicit function theorem entails that $c(t)$ is
differentiable in a neighborhood of $t=0$ and
$$ c'(0) = \frac{\e_{x_0}\etc{ \sum_0^{\tau(x_0)-1} f(S_i)}}{\e_{x_0}\etc{ \tau(x_0)}}\,.$$
Since $f$ is bounded, hence $\mu$-integrable, the ergodic theorem
implies
$$c'(0)=\frac{\dE_{x_0}\SBRA{ \sum_{i=0}^{\tau-1}
    f(S_i)}}{\dE_{x_0}\SBRA{\tau}}=\lim_{n \to \infty}
\frac{1}{n}\dE_{x_0}\SBRA{ \sum_{i=1}^{n} f(S_i)}=\langle \mu,f
\rangle \,,$$ \noindent with $\langle \mu,f \rangle= \sum_{x \in
\Si} f(x) \mu(x).$

\end{eproof}

\let\ep=\epsilon

\bigskip
We now prove that the constant $c(f)$ appearing in
Proposition~\ref{pro:large-devi-result} does not really depend on
the starting point $x_0$. Let  $(S_n)$ is a Markov chain taking
values in a countable set $\Sigma$. For any $x\in\Sigma$ define
$$ \kappa(x)=\sup\ens{\alpha>0 : \mathbb{E}_x\etc{e^{\alpha \tau(x)}} < +\infty}\,,\quad \text{with } \tau(x)=\inf\ens{n\ge 1 : S_n=x}\,.$$

Let $f:\Sigma\to \R$ be a bounded function. If $\norme{f}_\infty <
\kappa(x)  $ then the following limit exists
$$ c(x,f)= \lim \unsur{n} \log \mathbb{E}_x\etc{e^{A_n}}\,,\quad \text{with \,} A_n=\sum_{i=0}^{n-1} f(S_n)\,.$$
Different state points $x,y$ need to communicate to have the same
coefficient.
\begin{elem}
  If $(S_n)$ is irreducible recurrent and $\norme{f}_\infty <  \undemi \inf(\kappa(x),\kappa(y))$ then $c(x,f)=c(y,f)$.
\end{elem}
\begin{proof}
  Since  $(S_n)$ is irreducible recurrent, $\PP_x\etp{\tau(y)<+\infty}=1$ and there exists $p\ge 1$ such that  $\PP_x\etp{\tau(y)=p}>0$. Thanks to the strong Markov property,
$$ \mathbb{E}_x\etc{e^{A_n}}\ge  \mathbb{E}_x\etc{e^{A_n}\un{\tau(y)=p}} =
\mathbb{E}_x\etc{e^{A_p}\un{\tau(y)=p}\mathbb{E}_y\etc{e^{A_{n-p}}}}\,.$$
Let $\ep>0$. There exists $n_0$ such that for all $n\ge n_0$,
$\unsur{n} \mathbb{E}_y\etc{e^{A_{n}}} \ge c(y,f)-\epsilon$.

Therefore, if $n\ge n_0+p$, then
$$\mathbb{E}_x\etc{e^{A_n}}\ge e^{n (c(y,f)-\ep)} \mathbb{E}_x\etc{e^{A_p}\un{\tau(y)=p}}$$
and this yields
$$c(x,f)=\lim \unsur{n}  \log \mathbb{E}_x\etc{e^{A_n}} \ge c(y,f) -\ep\,.$$
Letting $\ep\to 0$ we get $c(x,f)\ge c(y,f)$. Substituting $x$ for
$y$, we obtain $c(x,f)=c(y,f)$.
\end{proof}

\begin{erem}\label{R:9}
  With the same argument we can prove that  $c(x,\beta,f)$ is the same for  the starting points $x$ ans $y$, as soon as $\l(\beta)+ 2 \norme{f}_\infty <  \inf(\kappa(x),\kappa(y))$.
\end{erem}

\section{Integration by parts formula for infinitely divisible laws}

Recall that the random variable $g$ has small exponential moments.
We assume now that it is infinitely divisible, and hence we have a
Levy Khinchine formula
\begin{equation}
  \label{eq:lkf}
   \l(\beta)=\log \gesp{e^{\beta g}} = c \beta + \frac{\sigma^2}{2} \beta^2 + \int \pi(du)\etp{ e^{\beta u} - 1 - \un{\valabs{u}\le 1} \beta u}\quad (\valabs{\beta}<\beta_0)\,,
\end{equation}

where $c\in\R$, $\sigma \ge 0$ are constants and $\pi$ is a
measure on $\R\backslash\ens{0}$ satisfying $\int \pi(du)\,
(1\wedge u^2) < +\infty$. \label{general_ipp}
\begin{elem} \label{lemme_ipp}
If $g$ satisfies \eqref{eq:lkf}, then for any bounded
differentiable $f$ with bounded derivative, one has the following
integration by parts formula:
\begin{equation} \label{ipp}
\gesp{gf(g)}=c \gesp{ f(g)} + \si^2 \gesp{ f'(g)} +
\int_{-\infty}^{+\infty}  \pi(du) u \SBRA{ \gesp{ f(g+u)} -
\un{\ABS{u}\le 1} \gesp{ f(g)}}
\end{equation}
\end{elem}
\begin{eproof}
As pointed out by  Nicolas Privault,   this Lemma can be seen as
an easy consequence of much more general integration by parts
formulas on the Poisson space (see Picard~\cite{MR98c:60055}). Let
us give a short proof here: it suffices  to prove the formula
(\ref{ipp}) for $f(x)=e^{i\te x}$, the extension to more general
functions following standard arguments. In that case,
$\gesp{gf(g)}=\gesp{ge^{i\te g}}= e^{\la(i\te)} \la'(i\te)$. Since
$\la'(i\te)=c+ \si^2 i\te + \int_{-\infty}^{+\infty}  \pi(du) u
\PAR{e^{i\te u}-\un{\ABS{u}\le 1}}$, we obtain:
\begin{eqnarray*}
e^{\la(i\te)} \la'(i\te)&=& c e^{\la(i\te)} + \si^2 i \te
e^{\la(i\te)} + \gesp{ e^{i\te g}} \int_{-\infty}^{+\infty}
\pi(du) u \PAR{e^{i\te
    u}-\un{\ABS{u}\le 1}} \cr
 &=& c \gesp{ f(g)} + \si^2 \gesp{ f'(g)} +
\int_{-\infty}^{+\infty}  \pi(du) u \SBRA{ \gesp{ f(g+u)} -
\un{\ABS{u}\le 1} \gesp{ f(g)}}.
\end{eqnarray*}
\end{eproof}

We shall now link the derivative of the free energy to
$\moynbis{L_n(S^1,S^2)}$, where  here and in the sequel, $L_n(S^1,
S^2) \df \sum_{i=1}^n 1_{(S^1_i= S^2_i)}$ denotes the global
correlation between the two independent configurations $S^1$ and
$S^2$ (under the same polymers measure $\langle \cdot
\rangle^{(n)}$). Recall that $I$ is an open interval, chosen as
big as possible, such that $0\in I \subset \ens{\beta : \l(\beta)<
+\infty}$.

\begin{eprop} \label{basee_sur_IPP}
If $g$ satisfies \eqref{eq:lkf}, then  there exists $c_1>0$,
depending on the law of $g$ and on $\be$, such that $\forall \be
\in I \cap (0, \infty)$
\begin{equation} \label{majoration-derivee}
 p'_n(\be)\le -\frac{c_1}{n} \gesp{ \moynbis{L_n(S^1,S^2)}}.
\end{equation}
\noindent In particular, for all $n\ge 1$, $\be \mapsto p_n(\be)$
is non increasing.
\\
Moreover if $2\be \in I \cap (0, \infty)$, there exists
$c_2(\beta)>0$ such that
\begin{equation} \label{minoration-derivee}
 p'_n(\be) \ge -\frac{c_2}{n} \gesp{\moynbis{L_n(S^1,S^2)}}.
\end{equation}

\end{eprop}

\begin{eproof}
In the sequel we write $\moybis{.}$ instead of $\moyn{.}$. The
first step is the following identity:
\begin{equation} \label{debut}
n p'_n(\be) = \gesp{\moybis{ \sum_{i=1}^n g(i,S_i) }} - n
\la'(\be)
 =  \sum_{i,x} \gesp{ \SBRA{ g(i,x) F_{i,x} (g(i,x))}} - n \la'(\be)  \, ,
\end{equation}
\noindent where we have set, for each $(i,x)$:
$$
F_{i,x}(u) \df \frac{ \dE \SBRA{ \un{S_i=x} \exp\PAR{\be
\sum_{(j,y)\neq (i,x)} g(j,y) \un{S_j=y} + \be u \un{S_i=x} }}}
{\dE \SBRA{ \exp\PAR{\be \sum_{(j,y)\neq (i,x)} g(j,y) \un{S_j=y}
+ \be u \un{S_i=x} }}}\, ,  \quad u\in \r.
$$
\noindent Since $F_{i,x}$ is a random function depending only on
 $(g(j,y),(j,y)\neq (i,x))$, it is independent of $g(i,x)$, so by Lemma
\ref{lemme_ipp}, one has for each fixed $(i,x)$:
\begin{align} \label{application_ipp}
\gesp{ g(i,x) F_{i,x} (g(i,x))} &= c \gesp{ F_{i,x} (g(i,x))}
+ \si^2 \gesp{ F'_{i,x} (g(i,x))}  \nonumber\\
&+\int_{-\infty}^{+\infty}  \pi(du) u \SBRA{ \gesp{
F_{i,x}(g(i,x)+u)} - \un{\ABS{u}\le 1} \gesp{ F_{i,x} (g(i,x))} }
\end{align}

\noindent Here one easily obtains that $F'_{i,x}(u)= \be
F_{i,x}(u) \SBRA{1-F_{i,x}(u)}$. In particular, one has
$$F'_{i,x}(g(i,x))=\be \moybis{\un{S_i=x}} \PAR{1-\moybis{\un{S_i=x}}}.$$
\noindent Moreover,
$$F_{i,x}(g(i,x)+u)=\frac{\moybis{\un{S_i=x}}e^{\be u}} {
\moybis{\un{S_i\neq x}}+ \moybis{\un{S_i=x}}e^{\be u}} \; ,$$
\noindent so that formula (\ref{application_ipp}) leads to:
\begin{eqnarray*}
\gesp{g(i,x) F_{i,x}(g(i,x))}&=& c \gesp{ \moybis{\un{S_i=x}}}  +
\si^2 \be \gesp{\moybis{\un{S_i=x}} \PAR{1-\moybis{\un{S_i=x}}}}
\cr &+& \int_{-\infty}^{+\infty}  \pi(du) u \SBRA{ \gesp{ \frac{
\moybis{\un{S_i=x}} e^{\be u} } { \moybis{\un{S_i\neq x}} +
\moybis{\un{S_i=x}} e^{\be u} } } - \un{\ABS{u}\le 1}
\gesp{\moybis{\un{S_i=x}}}  }
\end{eqnarray*}
\noindent Then, using that $\la'(\be)=c + \si^2 \be +
\int_{-\infty}^{+\infty} \pi(du) u (e^{\be u}- \un{\ABS{u}\le 1})$
and remembering that $ \sum_{i,x} \moybis{\un{S_i=x}} = n$,
equation (\ref{debut}) becomes
\begin{eqnarray*}
n p'_n(\be) &=& -\si^2 \be \gesp{\sum_{i,x} \moybis{\un{S_i=x}}^2}
- \sum_{i,x} \int_{-\infty}^{+\infty}  \pi(du) u \gesp{ \frac{
\moybis{\un{S_i=x}}^2 e^{\be u} (e^{\be u}-1) }
 {1 + \moybis{\un{S_i=x}} (e^{\be u}-1)}}.
\end{eqnarray*}
\noindent Now we prove  that
\begin{equation} \label{technique_hu}
\inf_{0 \le a \le 1} \int_{-\infty}^{+\infty}  \pi(du) u \frac{
e^{\be u} (e^{\be u}-1) } {1 + a (e^{\be u}-1)} >0
\end{equation}
as soon as $\pi(.) \neq 0$. On the one hand, if $supp(\pi) \cap
\dR_+ \neq \emptyset$, then for all $0 \le a \le 1$,
$$\int_{0}^{+\infty}  \pi(du) u
\frac{ e^{\be u} (e^{\be u}-1) } {1 + a (e^{\be u}-1)} \ge
\int_{0}^{+\infty}  \pi(du) u (e^{\be u}-1)   > 0 \, .$$ On the
other hand, if $supp(\pi) \cap \dR_- \neq \emptyset$, then for all
$0 \le a \le 1$,
$$\int_{-\infty}^0  \pi(du) \ABS{u}
\frac{ e^{\be u} (1-e^{\be u}) } {1 - a (1-e^{\be u})} \ge
\int_{-\infty}^0  \pi(du) e^{\be u} \ABS{u} (1-e^{\be u})   > 0 \,
.$$ \noindent In all cases, (\ref{technique_hu}) is true for all
$\pi(.) \neq 0$, so there exists $c_1> 0$ such that
$$n p'_n(\be) \le -c_1 \gesp{ \sum_{i,x} \moybis{\un{S_i=x}}^2} \, ,$$
\noindent $c_1$ being positive because $\pi(.) \neq 0$ or $\si>0$,
since the law of $X$ is non-degenerate. This leads the upper bound
(\ref{majoration-derivee}) thanks to the following identity:
\begin{equation} \label{reecriture_tpslocal}
\sum_{i,x}\moybis{\un{S_i=x}}^2=\moynbis{L_n(S^1,S^2)} \, .
\end{equation}
The lower bound (\ref{minoration-derivee}) can be deduced in the
same way, using that
\begin{eqnarray*}
\sup_{0 \le a \le 1} \int_{-\infty}^{+\infty}  \pi(du) u \frac{
e^{\be u} (e^{\be u}-1) } {1 + a (e^{\be u}-1)}
&\le& \int_{0}^{+\infty}  \pi(du) u  e^{\be u} (e^{\be u}-1)\\
&+& \int_{-\infty}^0  \pi(du) \ABS{u} (1-e^{\be u}) \\
&\stackrel{def}{=}& c' \, ,
\end{eqnarray*}
because $c'< \infty$ provided that $2\be \in I$.
\end{eproof}

\section{Strong disorder: Proof of Theorem \ref{countable_case}}

The part (a) of Theorem  \ref{countable_case} follows from
Proposition \ref{P:p4}, whereas the part (b) from Proposition
\ref{basee_sur_IPP}. To show the part (c), we make use of the
monotonicity of $\beta \to p_n(\beta)$, then  it suffices to prove
that $\limsup p_n(\beta)<0$ for $\beta >0$ small enough. Recall
that ${\bf E}(g)=\la'(0)=0$. Then for all $q=(q(i,x), i \ge 1, x
\in \Si) \in \dR^{\dN^*\times \Si}$, one has, using Jensen's
inequality:
\begin{align*}
p_n(\be)&=\frac{1}{n} \gesp{ \log \espxo{\PAR{
e^{\be \sum_{i,x} g(i,x)  \ind_{S_i=x} -n\la(\be) } }} - \beta \, \sum_{i,x} g(i,x) q(i,x)} \\
&=\frac{1}{n} \gesp{ \log \espxo{\PAR{
e^{\be \sum_{i,x} g(i,x) \PAR{ \ind_{S_i=x}-q(i,x)}-n\la(\be) } }}} & \\
&\le \frac{1}{n} \log \dE_{x_0} \PAR{
  e^{\sum_{i,x} \la\PAR{\be\SBRA{\ind_{S_i=x}-q(i,x)}}-n\la(\be) } }\, ,
\end{align*}

\noindent Let us choose $q(i,x)=\mu(x)$ (the invariant probability
measure of $S$) for all $(i,x) \in \dN^*\times \Si$ and let us fix
$\ep>0$. There exists $\be_{\ep}>0$ such that $\forall \, 0 < \be
< \be_{\ep}$, $\frac{1-\ep}{2} \be^2 \la''(0) \le \la(\be) \le
\frac{1+\ep}{2} \be^2 \la''(0)$. Thus, for $0 < \be < \be_{\ep}$,
\begin{eqnarray*}
p_n(\be) &\le& \frac{1}{n}\log \dE_{x_0} \PAR{e^{\frac{1+\ep}{2}
\be^2 \la''(0) \sum_{i,x}\PAR{ \ind_{S_i=x}-\mu(x)}^2 -
n\frac{1-\ep}{2} \be^2 \la''(0)}} \cr &\le&  \frac{1}{n}\log
\dE_{x_0} \PAR{e^{-(1+\ep) \be^2 \la''(0) \sum_{i=1}^n \mu(S_i)}}
+\ep \la''(0) \be^2 + \frac{1+\ep}{2} \be^2 \la''(0) \NRM{\mu}^2,
\end{eqnarray*}
\noindent with $\NRM{\mu}^2=\sum_{ x \in \Si} \mu^2(x)$.
\\
\noindent Applying  Lemma \ref{travail_Carmona} to $f(x)=-\mu(x)$,
one deduces the existence of $J=\SBRA{-t_0,t_0}$ and of $c_{\mu}$
defined on $J$ such that $ \forall \, t \in J, \; \frac{1}{n}\log
\dE_{x_0} \PAR{e^{t\sum_{i=1}^n \mu(S_i)}} \xrightarrow[n\to
\infty]{} c_{\mu}(t)$ and $\forall \, t \in J, \; c_{\mu}(t) \le
-(1-\ep)\ t \NRM{\mu}^2$, since $c'_{\mu}(0)=-\NRM{\mu}^2$. Hence
for $\be$ small enough one concludes that $\limsup_{n \to \infty}
p_n(\be) \le -\frac{1}{2} \be^2 \NRM{\mu}^2 \la''(0)$ and thus
$\limsup_{n \to \infty} p_n(\be) <0$ since $\la''(0)=\mbox{Var}
(g) >0$. Finally,  for any $0<\beta< \beta_0$, we deduce from the
property of concentration of measure (Lemma \ref{L:cm}, (ii)) that
almost surely, $Z_n(\beta)$ converges to $0$ exponentially fast.
\hfill $\Box$


\begin{thebibliography}{CSY03}

\bibitem{DeA98}
A. de Acosta  and  P. Ney.  \newblock Large deviation lower bounds
for arbitrary additive functionals of a Markov chain.
\newblock {\em Ann. Probab. } 26 (1998), no. 4, 1660--1682.


\bibitem{carmona-hu}
Ph. Carmona and Y. Hu.
\newblock On the partition function of a directed polymer in a {G}aussian
  random environment.
\newblock {\em Probab. Theory Related Fields}, 124(3):431--457, 2002.

\bibitem{comets}
F.~Comets, T.~Shiga, and N.~Yoshida.
\newblock Directed polymers in random environment: path localization and strong
  disorder.
\newblock {\em Bernoulli } (2003)   9(4):705--723, 2003.


\bibitem{comets-continu}
F.~Comets and N.~Yoshida.
\newblock Brownian directed polymers in random environment.
\newblock {\em preprint}, 2003.


\bibitem{DZ98}
A. Dembo  and  O. Zeitouni. {\it Large deviations techniques and
applications.} Second edition.  Springer-Verlag, New York, 1998.

\bibitem{guerra-toninelli}
F. Guerra and F.~L. Toninelli.
\newblock Quadratic replica coupling in the {S}herrington-{K}irkpatrick mean
  field spin glass model.
\newblock {\em J. Math. Phys.}, 43(7):3704--3716, 2002.

\bibitem{imbrie-spencer}
J.~Z. Imbrie and T.~Spencer.
\newblock Diffusion of directed polymers in a random environment.
\newblock {\em J. Statist. Phys.}, 52(3-4):609--626, 1988.

\bibitem{NN87}
P. Ney and E. Nummelin. Markov additive processes. I. Eigenvalue
properties and limit theorems. {\it Ann. Probab.}  15 (1987), no.
2, 561--592.

\bibitem{petermann}
M. Petermann.
\newblock Superdiffusivity of directed polymers in a random environment.
\newblock part of PHD thesis, 2000.

\bibitem{MR98c:60055}
J. Picard.
\newblock Formules de dualit\'e sur l'espace de {P}oisson.
\newblock {\em Ann. Inst. H. Poincar\'e Probab. Statist.}, 32(4):509--548,
  1996.



\bibitem{sinai}
Y.~Sinai.
\newblock A remark concerning random walks with random potentials.
\newblock {\em Fundamenta Mathematicae}, 147(2):173--180, 1995.

\end{thebibliography}

\end{document}